\UseRawInputEncoding
\documentclass[12pt]{article}

\usepackage{dsfont}
\usepackage{amsfonts,amsmath,amsthm}
\usepackage{algorithm, algorithmic}
\usepackage{color}
\usepackage{graphicx}
\usepackage{stmaryrd}        

\usepackage[colorlinks=true]{hyperref}
\hypersetup{urlcolor=blue, citecolor=blue,linkcolor=blue}

\usepackage[paper=a4paper,dvips,top=2cm,left=2cm,right=2cm,
    foot=1cm,bottom=4cm]{geometry}

\begin{document}

\title{von Neumann type trace inequality for dual quaternion matrices}
\author{ Chen Ling,\footnote{%
    Department of Mathematics, Hangzhou Dianzi University, Hangzhou, 310018, China
    ({\tt
macling@hdu.edu.cn}). This author's work was supported by Natural Science Foundation of China (No. 11971138).}
    \and
    Hongjin He,\thanks{School of Mathematics and Statistics, Ningbo University, Ningbo, 315211, China.
 ({\tt  hehongjin@nbu.edu.cn}). This author's work was supported by Natural Science Foundation of China (No. 11771113).}
\and
Liqun Qi\thanks{Department of Applied Mathematics, The Hong Kong Polytechnic University, Hung Hom,
    Kowloon, Hong Kong; Department of Mathematics, Hangzhou Dianzi University, Hangzhou, 310018, China ({\tt maqilq@polyu.edu.hk}).}
\and  Tingting Feng\thanks{%
    Department of Mathematics, Hangzhou Dianzi University, Hangzhou, 310018, China
    ({\tt
ttfeng@hdu.edu.cn}). This author's work was supported by Natural Science Foundation of China (No.12001146) and  the Zhejiang Provincial Natural Science Foundation of China (No. LQ21A010006).}
}
\date{\today}
\maketitle

\begin{abstract}
Dual quaternion matrices have important applications in multi-agent formation control. In this paper, we first address
the concept of spectral norm of dual quaternion matrices. Then, we introduce a von Neumann type trace inequality and a Hoffman-Wielandt type inequality for general dual quaternion matrices, where the latter characterizes a simultaneous perturbation bound on all singular values of a dual quaternion matrix. In particular, we also present two variants of the above two inequalities expressed by eigenvalues of dual quaternion Hermitian matrices. Our results are helpful for the further study of dual quaternion matrix theory, algorithmic design, and applications.

\medskip



\end{abstract}

\renewcommand{\Re}{\mathds{R}}
\newcommand{\rank}{\mathrm{rank}}
\renewcommand{\span}{\mathrm{span}}
\newcommand{\X}{\mathcal{X}}
\newcommand{\A}{\mathcal{A}}
\newcommand{\I}{\mathcal{I}}
\newcommand{\B}{\mathcal{B}}
\newcommand{\C}{\mathcal{C}}
\newcommand{\OO}{\mathcal{O}}
\newcommand{\e}{\mathbf{e}}
\newcommand{\0}{\mathbf{0}}
\newcommand{\dd}{\mathbf{d}}
\newcommand{\ee}{\mathbf{e}}
\newcommand{\ii}{\mathbf{i}}
\newcommand{\jj}{\mathbf{j}}
\newcommand{\kk}{\mathbf{k}}
\newcommand{\va}{\mathbf{a}}
\newcommand{\vb}{\mathbf{b}}
\newcommand{\vc}{\mathbf{c}}
\newcommand{\vg}{\mathbf{g}}
\newcommand{\vr}{\mathbf{r}}
\newcommand{\vt}{\rm{vec}}
\newcommand{\vx}{\mathbf{x}}
\newcommand{\vy}{\mathbf{y}}
\newcommand{\vu}{\mathbf{u}}
\newcommand{\vv}{\mathbf{v}}
\newcommand{\y}{\mathbf{y}}
\newcommand{\vz}{\mathbf{z}}
\newcommand{\T}{\top}

\newtheorem{Thm}{Theorem}[section]
\newtheorem{Def}[Thm]{Definition}
\newtheorem{Ass}[Thm]{Assumption}
\newtheorem{Lem}[Thm]{Lemma}
\newtheorem{Prop}[Thm]{Proposition}
\newtheorem{Cor}[Thm]{Corollary}
\newtheorem{example}[Thm]{Example}
\newtheorem{ReK}[Thm]{Remark}
\newtheorem{PropT}[Thm]{Property}

\section{Introduction}\label{Introd}
Dual quaternions, which serve as the core part of Clifford algebra or geometric algebra,  were originally introduced by Clifford \cite{Cl73} in 1873. Mathematically, the dual quaternions are an $8$-dimensional real algebra isomorphic to the tensor product of the quaternions and the dual numbers. In mechanics, the dual quaternions give a neat and succint way to encapsulate both translations and rotaions into a unified representation for rigid transformations in three dimensions. With the rapid developments of artificial intelligence, dual quaternions have been successfully applied to the areas of automatic control, computer vision and bioengineering, such as robotics control \cite{Da99,MXXLY18,Tho14,WYL12,WZ14} and 3D computer graphics \cite{BLH19,CKJC16,TRA11}, neuroscience \cite{LLB13} and biomechanics \cite{PV10}, to name just a few.
More interestingly, we notice that dual quaternion matrices have important application in multi-agent formation control \cite{QWL22}. In 2011, Wang \cite{Wa11} raised this issue.   In an unpublished manuscript, Wang, Yu and Zheng \cite{WYZ} introduced three classes of dual quaternion matrices for studying multi-agent formation control, namely relative configuration adjacency matrices, logarithm adjacency matrices and relative twist adjacency matrices. Recently, Qi and Luo \cite{QL21} showed that dual quaternion Hermitian matrices have very nice spectral properties. They showed that an $n \times n$ dual quaternion Hermitian matrix has exactly $n$ right eigenvalues, which are all dual numbers and are also the left eigenvalues of this Hermitian matrix. Thus, we may simply call them eigenvalues of that Hermitian matrix. This Hermitian matrix is positive semi-definite or positive definite if and only if these $n$ dual number eigenvalues are nonnegative or positive in the sense of \cite{QLY21}, respectively. Moreover, the minimax principle and generalized inverses of dual quaternion matrices were studied in \cite{LQY22}, and the singular values of dual quaternion matrices and their low-rank approximations also were studied \cite{LHQ22}. Then, Qi, Wang and Luo \cite{QWL22} showed that the relative configuration adjacency matrix and the logarithm adjacency matrix are dual quaternion Hermitian matrices. These dual quaternion matrices have important applications in multi-agent formation control, and the related spectral and positive semi-definite properties pave the way for us to study stability issues of the multi-agent formation control problem. The work in \cite{QWL22} builds a bridge for the research on dual quaternion matrices and multi-agent formation control.    Furthermore, Cui and Qi \cite{CQ23} proposed a power method to computing eigenvalues of a dual quaternion Hermitian matrix, and applied it to the simultaneous location and mapping problem.


As we know, the well-known von Neumann trace inequality \cite{VN37} (see also \cite{Mir75}), which bounds the inner product of two matrices via the inner product of their singular value vectors, is the key inequality for the analysis of spectral functions, and plays a pivot role in the developments of low-rank matrix approximation theory and low-rank optimization, e.g., see \cite{HJ13}. Since quaternion matrices have wide applications in engineering \cite{CQZ22,CXZ20,Gri17,LM04}, the related spectral theory has been received considerable attention in recent years \cite{WLZZ18,XM15,Zh97}. Particularly, to derive quaternionic proximity operators for trace-norm regularized optimization problems arising from audio separation, Chan and Yang \cite{CY16} first proved that the von Neumann trace inequality still holds for quaternion matrices. As a combination of dual numbers and quaternions, dual quaternion matrices have been applied to multi-agent formation control \cite{QWL22}. In theoretical aspects, however, the existence of zero divisors, i.e., infinitesimal dual quaternion numbers makes analysis on dual quaternion matrices difficult. It is unknown whether the von Neumann trace inequality still holds for dual quaternion matrices. To answer such question, we first introduce the spectral norm for dual quaternion matrices in this paper. Then, we present a von Neumann type trace inequality for dual quaternion matrices, which then paves the way to establish a Hoffman-Wielandt type inequality characterizing a simultaneous perturbation bound for all singular values of a general dual quaternion matrix. It is worth pointing out that, when the dual quaternion reduces the quaternion, the von Neumann inequality obtained in this paper is exactly the one presented in \cite{CY16}, but our proof method is completely  different from the way used in  \cite{CY16}, even in the case of quaternions. Moreover, by considering the application of dual quaternion Hermitian matrices in multi-agent formation control, we also discuss the above two inequalities for dual quaternion Hermitian matrices. We believe that our results will enrich the theory of dual quaternion matrices, and they will be of benefit to further study of dual quaternion matrices, algorithmic design, and applications.

This paper is divided into five sections. In Section \ref{Prelim}, we present some preliminaries on dual numbers, quaternions, dual quaternions and dual quaternion algebra. In Section \ref{Spec-Norm}, we introduce the concept of the spectral norm for dual quaternion matrices, which is exactly the largest singular value of the involved dual quaternion matrix. In Section \ref{VN-TraceIn}, we present a von Neumann type trace inequality for dual quaternion matrices. In Section \ref{sec_HWineq}, we consider the extension of the well-known Hoffman-Wielandt inequality to dual quaternionic versions. Finally, some concluding remarks are drawn in Section \ref{Conclusion}.

\section{Preliminaries}\label{Prelim}

\subsection{Dual numbers}\label{Du-number}
Let $\mathbb{R}$ and $\widehat{\mathbb R}$ denote the field of the real numbers and the set of the dual numbers, respectively. A dual number $q\in \widehat{\mathbb R}$ has the form $q = q_{\sf st} + q_{\sf in}\epsilon$, where $q_{\sf st}\in \mathbb{R}$ is called the real part or the standard part; $q_{\sf in}\in \mathbb{R}$ represents the dual part or the infinitesimal part of $q$; and $\epsilon$ is the infinitesimal unit satisfying $\epsilon^2 = 0$.  Particularly, if the standard part $q_{\sf st}$ of $q$ is nonzero, i.e., $q_{\sf st} \not = 0$, we say that $q$ is appreciable; otherwise, we say that $q$ is infinitesimal. Note that the infinitesimal unit $\epsilon$ is commutative in multiplication with real numbers, complex numbers and
 quaternion numbers (see Section \ref{Qua_number}). The dual numbers form a commutative algebra of dimension two over the reals.

Now, we recall the recently introduced total order for dual numbers in \cite{QLY21}. Given two dual numbers $p = p_{\sf st} + p_{\sf in}\epsilon, q = q_{\sf st} + q_{\sf in}\epsilon \in \widehat{\mathbb R}$ with $ p_{\sf st}, p_{\sf in},q_{\sf st}, q_{\sf in}$ being real numbers, we say that $q < p$ if either $q_{\sf st} < p_{\sf st}$, or $q_{\sf st} = p_{\sf st}$ and $q_{\sf in} < p_{\sf in}$, and we say that $q = p$ if $q_{\sf st} = p_{\sf st}$ and $q_{\sf in} = p_{\sf in}$. Consequently, if $q > 0$, we say that $q$ is a positive dual number; and if $q \ge 0$, we say that $q$ is a nonnegative dual number. In what follows, we denote the set of nonnegative and positive dual numbers by $\widehat{\mathbb R}_+$ and $\widehat{\mathbb R}_{++}$, respectively. For given $p = p_{\sf st} + p_{\sf in}\epsilon, q = q_{\sf st} + q_{\sf in}\epsilon \in \widehat{\mathbb R}$, 
we have
\begin{equation} \label{e1}
p + q =p_{\sf st}+q_{\sf st} +(p_{\sf in}+q_{\sf in})\epsilon,~~~~pq = p_{\sf st}q_{\sf st} +(p_{\sf st}q_{\sf in}+p_{\sf in} q_{\sf st})\epsilon.\end{equation}
Following the definition in  \cite{QLY21}, the absolute value of $q=q_{\sf st}+q_{\sf in}\epsilon\in \widehat{\mathbb R}$ is defined by
\begin{equation}\label{e5}
|q| = \left\{ \begin{array}{ll}|q_{\sf st}| + {\rm sgn}(q_{\sf st})q_{\sf in}\epsilon, & {\rm if~}  q_{\sf st} \not = 0, \\
|q_{\sf in}|\epsilon, & {\rm otherwise},  \end{array}  \right.
\end{equation}
where `${\rm sgn}(\cdot)$' represents the sign function, that is, for any $u \in \mathbb R$,
$${\rm sgn}(u) =  \left\{ \begin{aligned} 1, & \ {\rm if}\  u > 0, \\ 0, &   \ {\rm if}\  u = 0, \\
-1, &   \ {\rm if}\  u < 0.  \end{aligned}  \right. $$
For a given dual number $q=q_{\sf st}+q_{\sf in}\epsilon$, if $q$ is appreciable, then $q$ is nonsingular and $q^{-1} = q_{\sf st}^{-1} - q_{\sf st}^{-1}q_{\sf in} q_{\sf st}^{-1}\epsilon$. If $q$ is infinitesimal, then $q$ is not nonsingular. If $q$ is nonnegative and appreciable, then the square root of $q$ is still a nonnegative dual number. If $q$ is positive and appreciable, we have
\begin{equation} \label{e4}
\sqrt{q} = \sqrt{q_{\sf st}} + {q_{\sf in} \over 2\sqrt{q_{\sf st}}}\epsilon.
\end{equation}
In particular, we have $\sqrt{q} = 0$ when $q=0$.

Below, we recall a proposition introduced in \cite{LQY22}.
\begin{Prop}\label{P6.5}
Let $p,q\in \widehat{\mathbb R}$. Then, we have the following conclusions.
\begin{itemize}
\item[{\rm(a).}] If $p, q \in \widehat{\mathbb R}_{+}$, then $pq \in \widehat{\mathbb R}_{+}$.
\item[{\rm(b).}] If $p, q \in \widehat{\mathbb R}_{++}$ and at least one of them is appreciable, then $pq \in \widehat{\mathbb R}_{++}$.
\item[{\rm(c).}] If $p\geq 0$, then $|p|=p$; otherwise, $|p|>p$.
\item[{\rm(d).}] If $p$ is appreciable, then $|p|=\sqrt{p^2}$.
\item[{\rm(e).}] If $p, q \in \widehat{\mathbb R}_{++}$ and are both appreciable, then $\sqrt{pq}=\sqrt{p}\sqrt{q}$.
\item[{\rm(f).}] If $q\in \widehat{\mathbb R}_{++}$ and is appreciable, then $p-q \in \widehat{\mathbb R}_{+}$ implies $\sqrt{p}-\sqrt{q} \in \widehat{\mathbb R}_{+}$.
\end{itemize}
\end{Prop}

\subsection{Quaternions}\label{Qua_number}
Denote by $\mathbb Q$ the four-dimensional vector space of the quaternions over $\mathbb{R}$, with an ordered basis,
denoted by $\ee, \ii, \jj$ and $\kk$. A quaternion $q\in \mathbb{Q}$ has the form
$q = q_0\ee + q_1\ii + q_2\jj + q_3\kk$, where $q_0, q_1, q_2$ and $q_3$ are real numbers, $\ii, \jj$ and $\kk$ are three imaginary units of quaternions, which satisfy
$$\ii^2 = \jj^2 = \kk^2 =\ii\jj\kk = -1,\quad \ii\jj = -\jj\ii = \kk, \quad \jj\kk = - \kk\jj = \ii, \quad  \kk\ii = -\ii\kk = \jj.$$
For notational simplicity, we will omit the real unit $\ee$ and denote $q \in \mathbb Q$ as $q = q_0 + q_1\ii + q_2\jj + q_3\kk$ in the subsequent discussions.

For a quaternion $q = q_0 + q_1\ii + q_2\jj + q_3\kk$, we call ${\rm Re}(q):= q_0$ and ${\rm Im}(q):= q_1\ii + q_2\jj +q_3\kk$ the real and imaginary parts, respectively; the conjugate of $q$ is given by $\bar q := q_0 -q_1\ii -q_2\jj -q_3\kk$ and the norm of $q$ is defined by $|q|:=\sqrt{\bar qq}=\sqrt{q_0^2+q_1^2+q_2^2+q_3^2}$.
In particular, a quaternion is called imaginary if its real part is zero. Given two quaternions $p = p_0 + p_1\ii + p_2\jj + p_3\kk$ and $q = q_0 + q_1\ii + q_2\jj + q_3\kk$, it is easy to verify that
\begin{equation}\label{pqmut}
p\bar q+q\bar p=\bar pq +\bar qp = 2(p_0q_0+p_1q_1+p_2q_2+p_3q_3).
\end{equation}
However, we shall notice that the multiplication of quaternions satisfies the distribution law, but is noncommutative. In fact, $\mathbb{Q}$ is an associative but non-commutative algebra of four rank over $\mathbb{R}$, called quaternion skew-field \cite{WLZZ18}.

Throughout this paper, we denote by $\mathbb{Q}^{m\times n}$  the collection of all $m\times n$ matrices with quaternion entries. Specially, $\mathbb{Q}^{m\times 1}$ is abbreviated as $\mathbb{Q}^m$, which is the collection of quaternion column vectors with $m$ components.

\subsection{Dual quaternions and dual quaternion algebra}
\subsubsection{Dual quaternions}
Dual quaternion is a composite concept, which is the combination of dual numbers and quaternions. Specifically, a dual quaternion $q$ has the form $q = q_{\sf st} + q_{\sf in}\epsilon$,
where $q_{\sf st}, q_{\sf in} \in \mathbb {Q}$ are the standard part and the infinitesimal part of $q$, respectively. Throughout, we denote by $\widehat{\mathbb Q}$ the set of dual quaternions. Recalling the definitions introduced in \cite{BK20, CKJC16, Ke12}, for any two dual quaterions $p=p_{\sf st} + p_{\sf in}\epsilon$ and $q=q_{\sf st} + q_{\sf in}\epsilon$, the addition and multiplication between them are defined by $$p+q=(p_{\sf st}+q_{\sf st})+(p_{\sf in}+q_{\sf in})\epsilon$$
and
$$pq=p_{\sf st}q_{\sf st} + (p_{\sf in}q_{\sf st}+p_{\sf st}q_{\sf in})\epsilon,$$
respectively. It is easy to see that $\widehat{\mathbb Q}$ is a ring with respect to the two binary algebraic operations defined above. The conjugate of $q$ is $\bar q = \bar q_{\sf st} + \bar q_{\sf in}\epsilon$. It is easy to see that $q\bar q=\bar q q$ for any $q\in \widehat{\mathbb Q}$ \cite{QLY21}. Similar to dual numbers, if $q_{\sf st} \not = 0$, then we say that $q$ is appreciable, otherwise, we say that $q$ is infinitesimal. We can derive that $q$ is invertible if and only if  $q$ is appreciable. In this case, we have $q^{-1} = q_{\sf st}^{-1} - q_{\sf st}^{-1}q_{\sf in} q_{\sf st}^{-1} \epsilon$. The magnitude of $q\in \widehat{\mathbb Q}$ is defined as
\begin{equation} \label{e7}
\displaystyle|q| := \left\{ \begin{array}{ll} |q_{\sf st}| +\displaystyle {(q_{\sf st}\bar q_{\sf in}+q_{\sf in} \bar q_{\sf st}) \over 2|q_{\sf st}|}\epsilon, & \ {\rm if}\  q_{\sf st} \not = 0, \\
|q_{\sf in}|\epsilon, &  \ {\rm otherwise},
\end{array} \right.
\end{equation}
which is a dual number by (\ref{pqmut}). Notice that such a definition immediately reduces to the absolute function \eqref{e5} when $q$ is a dual number, i.e., $q \in \widehat{\mathbb R}$, and it is exactly the magnitude of a quaternion when $q \in \mathbb Q$.

\subsubsection{Dual quaternion algebra}
Denote by $\widehat{\mathbb Q}^{m}$ the set of all dual quaternion vectors with $m$ components. For any two $m$-dimensional dual quaternion vectors ${\bf x}=(x_1,x_2,\ldots,x_m)^\top$, ${\bf y}=(y_1,y_2,\ldots,y_m)^\top$ and a dual quaternion $\alpha\in \widehat{\mathbb Q}$, we define
$${\bf x}+{\bf y}=(x_1+y_1,x_2+y_2,\ldots,x_m+y_m)^\top\quad \text{and}\quad {\bf x}\alpha=(x_1\alpha,x_2\alpha,\ldots,x_m\alpha)^\top,$$
which is called the right multiplication of ${\bf x}\in \widehat{\mathbb Q}^m$ and $\alpha\in \widehat{\mathbb Q}$. It is easy to verify that $\widehat{\mathbb Q}^m$ is an $m$-dimensional vector space over $\widehat{\mathbb Q}$, with respect to the addition and right multiplication defined above.
\begin{Def}[\cite{LQY22}]\label{Def-DQVecIn}
Let $\Xi:=\{{\bf u}^{(1)},{\bf u}^{(2)},\ldots,{\bf u}^{(s)}\}\subset \widehat{\mathbb Q}^m$. We say that $\Xi$ is right linearly independent, if for any $\alpha_1,\alpha_2,\ldots,\alpha_s\in \widehat{\mathbb Q}$,
$$
{\bf u}^{(1)}\alpha_1+{\bf u}^{(2)}\alpha_2+\ldots+{\bf u}^{(s)}\alpha_s={\bf 0}~~~\Rightarrow~~~\alpha_1=\alpha_2=\ldots=\alpha_s=0.
$$
\end{Def}
As a result of Definition \ref{Def-DQVecIn}, we can see that, if $\Xi$ is right linearly independent, then ${\bf u}^{(i)}$ is appreciable for every $i=1,2,\ldots,s$.  
For given ${\bf u}=(u_1,u_2,\ldots,u_m)^\top$ and ${\bf v}=(v_1,v_2,\ldots,v_m)^\top$ in $\widehat{\mathbb Q}^m$, denote by $\langle {\bf u},{\bf v}\rangle$ the dual quaternion-valued inner product, i.e., $\langle {\bf u},{\bf v}\rangle=\sum_{i=1}^m\bar v_iu_i$. It is easy to see that $\langle {\bf u}, {\bf v} \alpha+{\bf w} \beta\rangle= \bar\alpha\langle {\bf u}, {\bf v}\rangle +\bar\beta\langle{\bf u},{\bf w}\rangle$ and $\langle {\bf u},{\bf v}\rangle=\overline{\langle {\bf v},{\bf u}\rangle}$ for any ${\bf u},{\bf v}, {\bf w}\in \widehat{\mathbb Q}^m$ and $\alpha,\beta\in \widehat{\mathbb Q}$. If ${\bf u}\in \widehat{\mathbb Q}^m$ is appreciable, then $\langle {\bf u},{\bf u}\rangle$ is an appreciable positive dual number.

\begin{Def}
\label{Def-Orthog}
Let ${\bf u},{\bf v}\in \widehat{\mathbb Q}^m$ be appreciable. We say that ${\bf u},{\bf v}$  are orthogonal if $\langle {\bf u},{\bf v}\rangle=0$. An $n$-tuple $\{{\bf u}^{(1)},{\bf u}^{(2)},\ldots, {\bf u}^{(s)}\}\subset \widehat{\mathbb Q}^m$, where all ${\bf u}^{(1)},{\bf u}^{(2)},\ldots, {\bf u}^{(s)}$  are appreciable, is said to be orthogonal if $\langle {\bf u}^{(i)},{\bf u}^{(j)}\rangle=0$ for $i\neq j$, and orthonormal if it is orthogonal and $\langle {\bf u}^{(i)},{\bf u}^{(i)}\rangle=1$ for $i=1,2,\ldots,s$.
\end{Def}

Denote by $\widehat{\mathbb Q}^{m\times n}$ the set of $m\times n$ dual quaternion matrices. Then $A=(a_{ij})\in \widehat{\mathbb Q}^{m\times n}$ can
be written as $A = A_{\sf st} + A_{\sf in}\epsilon$, where $A_{\sf st}, A_{\sf in}\in \mathbb{Q}^{m\times n}$ are the standard part and the infinitesimal part of $A$, respectively.  If $A_{\sf st}$ is nonzero, i.e., $A_{\sf st} \not = O$, we say that $A$ is appreciable, otherwise, we say that $A$ is infinitesimal. For given $A\in \widehat{\mathbb Q}^{m\times n}$, the transpose of $A$ is denoted as $A^\top = (a_{ji})$, the conjugate of $A$ is denoted as ${\bar A} = (\bar a_{ij})$, and the conjugate transpose of $A$ is denoted as $A^*=(\bar a_{ji})=\bar A^\top$. It is obvious that $A^\top = A_{\sf st}^\top + A_{\sf in}^\top\epsilon$, $\bar A = \bar A_{\sf st}+\bar A_{\sf in}\epsilon$ and $A^* = A_{\sf st}^* + A_{\sf in}^*\epsilon$.  In this paper, a square matrix $A\in \widehat{\mathbb Q}^{m\times m}$ is called nonsingular (invertible) if $AB = BA = I_m$ for some $B\in \widehat{\mathbb Q}^{m\times m}$. In that situation, we denote $A^{-1} = B$. Moreover, a square matrix $A\in \widehat{\mathbb Q}^{m\times m}$ is called normal if $AA^*=A^*A$, Hermitian if $A^*=A$, and unitary if $A$ is nonsingular and $A^{-1}=A^*$. We have $(AB)^{-1}=B^{-1}A^{-1}$ if $A$ and $B$ are nonsingular, and $(A^*)^{-1}=(A^{-1})^*$ if $A$ is nonsingular. We say that $A\in \widehat{\mathbb Q}^{m\times m}$ is unitary, if $A$ satisfies $A^*A=I_m$. It is obvious that $A$ is unitary, if and only if the set consisting of column (row) vectors form an orthonormal basis of $\widehat{\mathbb Q}^m$, i.e., it is orthonormal and any vector in $\widehat{\mathbb Q}^m$ can be written as a right linear combination of this set. Similarly, we say that $A\in \widehat{\mathbb Q}^{m\times s} (s\leq m)$ is partially unitary, if $A$ satisfies $A^*A=I_s$. From Definition \ref{Def-DQVecIn}, it is easy to see that, the right linear independence of the involved vector set $\Xi$ is essentially that  $A{\bf x}={\bf 0}$ has a unique zero solution in $\widehat{\mathbb Q}^n$,  where $A=[{\bf u}^{(1)},\ldots,{\bf u}^{(s)}]\in \widehat{\mathbb Q}^{m\times s}$. For given $A=(a_{ij})\in \hat{\mathbb{Q}}^{m\times m}$, the trace of $A$, named ${\rm trace}(A)$, is defined as
\begin{equation}\label{TraceDef}
{\rm trace}(A)=\sum_{i=1}^ma_{ii},
\end{equation}
which is a dual quaternion.

As defined in \cite{QL21}, for a dual quaternion matrix $A\in \widehat{\mathbb Q}^{m\times m}$, if there exist a $\lambda\in \widehat{\mathbb Q}$ and an appreciable ${\bf x}\in \widehat{\mathbb Q}^m$ such that
\begin{equation}\label{RightEig}
A{\bf x} = {\bf x}\lambda,
\end{equation}
then we say that $\lambda$ is a right eigenvalue of $A$, with ${\bf x}$ as an associated right eigenvector. If $\lambda$ is a dual number, then we have
\begin{equation}\label{LeftEig}
A{\bf x} = \lambda{\bf x},
\end{equation}
i.e., $\lambda$ is also a left eigenvalue of $A$.   In this case, $\lambda$ is simply called an eigenvalue of $A$, and $\bf x$ an associated eigenvector.    In particular, it was shown in \cite{QL21} that an $m \times m$ dual quaternion Hermitian matrix has exactly $m$ dual number eigenvalues. 

For given $A = A_{\sf st} + A_{\sf in} \epsilon = (a_{ij}) \in {\widehat{\mathbb Q}}^{m \times n}$, the Frobenius norm of $A$, which is a dual number, is defined by
\begin{equation}\label{FNorm-DQM}
\|A \|_F = \left\{\begin{aligned}\sqrt{\sum_{i=1}^m \sum_{j=1}^n |a_{ij}|^2}, & \quad\ {\rm if}\ A_{\sf st} \not = O, \\
\|A_{\sf in}\|_F\epsilon,\quad & \quad \ {\rm otherwise.} \end{aligned}\right.
\end{equation}
Clearly, the Frobenius norm of a matrix is actually the $\ell_2$-norm of the vectorization of that matrix. Most recently, Ling et al. \cite{LQY22} proved a dual quaternion version of Cauchy-Schwarz inequality, which can be stated as follows.
\begin{Prop}[Cauchy-Schwarz inequality on $\widehat{\mathbb Q}^m$]\label{Ch-SW-In}
	For any ${\bf u},{\bf v}\in \widehat{\mathbb Q}^m$, it holds that $$\|{\bf u}\|_2\|{\bf v}\|_2-|\langle{\bf u},{\bf v}\rangle|\in \widehat{\mathbb R}_+,$$
	that is, $|\langle{\bf u},{\bf v}\rangle|\leq\|{\bf u}\|_2\|{\bf v}\|_2$.
\end{Prop}
As a consequence of (\ref{FNorm-DQM}) and Proposition \ref{Ch-SW-In}, for given $A \in {\widehat{\mathbb Q}}^{m \times n}$ and $\vx \in {\widehat{\mathbb Q}}^n$, regardless of whether $A\vx$ is appreciable or not,
it can be verified that
$\|A\vx\|_2 \le \|A\|_F \|\vx\|_2.$

\begin{Prop} \label{pp3.2}
Suppose that $U \in {\widehat{\mathbb Q}}^{m \times n}$ is partially unitary, and $\vx \in {\widehat{\mathbb Q}}^n$.   Then
\begin{equation} \label{eq5}
\|U\vx\|_2 = \|\vx\|_2.
\end{equation}
\end{Prop}
\begin{proof}   Suppose that $\vx = \vx_{\sf st}+\vx_{\sf in} \epsilon$ is appreciable. 
  It follows from (\ref{FNorm-DQM}) that
$$\|\vx\|_2^2 = \sum_{i=1}^n |x_i|^2 = \vx^*\vx.$$
On the other hand, let $U = U_{\sf st}+U_{\sf in}\epsilon$.  A direct calculation leads to $U_{\sf st}^* U_{\sf st}=I_n$.  
Then, the standard part of $U\vx$ is $U_{\sf st}\vx_{\sf st} \not = \0$, i.e., $U\vx$ is also appreciable. Consequently, by (\ref{FNorm-DQM}) again, we have
$$\|U\vx\|_2^2 = (U\vx)^*(U\vx) = \vx^*U^*U\vx = \vx^*\vx.$$
Hence, $\|U\vx\|_2 = \|\vx\|_2$ in this case.

Now,  we assume that $\vx$ is infinitesimal, i.e., $\vx = \vx_{\sf in}\epsilon$. Then $U\vx = U_{\sf st}\vx_{\sf in}\epsilon$, which means that $U\vx $ is also infinitesimal.  We have $\|U_{\sf st}\vx_{\sf in}\|_2 = \|\vx_{\sf in}\|_2$ since $U_{\sf st}^* U_{\sf st}=I_n$. Then by (\ref{FNorm-DQM}), we still have
$\|U\vx\|_2 = \|\vx\|_2$ in this case.
\end{proof}

\begin{Prop}[\cite{QLY21}]\label{p6.3}
For any $\vu =\vu_{\sf st} + \vu_{\sf in}\epsilon\in{\widehat{\mathbb Q}}^m$ with $\vu_{\sf st}\neq \0$, it holds that
\begin{equation}\label{e12-1}
\|\vu\|_2=\|\vu_{\sf st}\|_2+\displaystyle\frac{\langle\vu_{\sf st},\vu_{\sf in}\rangle+\langle\vu_{\sf in},\vu_{\sf st}\rangle}{2\|\vu_{\sf st}\|_2}\epsilon.
\end{equation}
\end{Prop}

\section{Spectral norm of dual quaternion matrices}\label{Spec-Norm}
We begin this section with recalling the following two theorems for dual quaternion matrices.

\begin{Thm}[\cite{QL21}]\label{HUDec}
	 Let $A=A_{\sf st}+A_{\sf in}\epsilon\in \widehat{\mathbb Q}^{m\times m}$ be a Hermitian matrix. Then, there exists a unitary matrix $U\in \widehat{\mathbb Q}^{m\times m}$ and a diagonal matrix $\Sigma\in \widehat{\mathbb Q}^{m\times m}$ such that $A=U\Sigma U^*$, where
\begin{equation}\label{Sigmn}\Sigma := {\rm diag} (\lambda_1(A),\lambda_2(A),\ldots,\lambda_m(A)),
\end{equation}
where $\lambda_1(A)\geq\lambda_2(A)\geq\ldots\geq\lambda_m(A)$ are dual numbers. Counting possible multiplicities $\lambda_{i,j}$, the form $\Sigma$ is unique.
\end{Thm}

It is obvious that $\lambda_i(A)$ is the $i$th largest eigenvalue of $A$, with ${\bf u}^{(i)}$ as an associated eigenvector, where ${\bf u}^{(i)}$ is the $i$th column in $U$. When $A$ is Hermitian, since $\lambda_1(A), \lambda_2(A), \ldots,\lambda_m(A)$ are dual numbers, from (\ref{Sigmn}) and (\ref{pqmut}), it is easy to see that
$$
a_{ii}=\sum_{j=1}^m\lambda_j(A)u_{ij}\bar u_{ij}=\sum_{j=1}^m\lambda_j(A)\bar u_{ij}u_{ij},
$$
which implies that
\begin{equation}\label{sumEigV}
\sum_{i=1}^ma_{ii}=\sum_{j=1}^m\lambda_j(A)\sum_{i=1}^m\bar u_{ij}u_{ij}=\sum_{j=1}^m\lambda_j(A)\|{\bf u}^{(j)}\|_2^2=\sum_{j=1}^m\lambda_j(A),
\end{equation}
since $\|{\bf u}^{(j)}\|_2=1$. Hence, if $A$ is Hermitian, we have ${\rm trace}(A)={\rm trace}(UAU^*)$ for any unitary matrix $U\in \widehat{\mathbb{Q}}^{m\times m}$.

\begin{Thm}[\cite{QL21}]\label{SVD-DQM}
For a given $A\in \widehat{\mathbb Q}^{m\times n}$, there exist two dual quaternion unitary matrices $U\in \widehat{\mathbb Q}^{m\times m}$ and $V\in \widehat{\mathbb Q}^{n\times n}$, such that
\begin{equation}\label{SVDDQMEQ}
A=U\left[\begin{array}{cc}\Sigma_t&O\\
O&O
\end{array} \right]_{m\times n}V^*,
\end{equation}
where $\Sigma_t\in \widehat{\mathbb R}^{t\times t}$ is a diagonal matrix, taking the form
$\Sigma_t={\rm diag} (\sigma_1(A),\ldots, \sigma_r(A),\ldots,\sigma_t(A))$, $r \leq t\leq s:={\min}\{m, n\}$, $\sigma_1(A)\geq \sigma_2(A)\geq\ldots\geq\sigma_r(A)$ are positive appreciable dual numbers, and $\sigma_{r+1}(A)\geq \sigma_{r+2}(A)\geq\ldots\geq\sigma_t(A)$ are positive infinitesimal dual numbers. Counting possible multiplicities of the diagonal entries, the form $\Sigma_t$ is unique.
\end{Thm}

We call form \eqref{SVDDQMEQ} the singular value decomposition (SVD) for dual quaternion matrix $A$, and call $\sigma_1(A),\ldots, \sigma_r(A),\ldots,\sigma_t(A)$ and possibly $\sigma_{t+1}(A)=\ldots=\sigma_s(A)=0$ (if $t< s)$ the singular values of
$A$, where $t$ and $s$ correspond to the rank and the appreciable rank of $A$, respectively.

For given $A\in {\widehat{\mathbb Q}}^{m \times n}$, the spectral norm $\|A\|_2$ is defined by
\begin{equation}\label{SNorm-DQM}
\|A \|_2 = \max_{{\bf x}\in \widehat{\mathbb Q}^n,~\|{\bf x}\|_2=1}\|A{\bf x}\|_2.
\end{equation}
Notice that $\|A\|_2$ is induced by the $\ell_2$-norm on dual quaternion vector spaces and hence is a matrix norm. In addition, by (\ref{SNorm-DQM}), we have
\begin{equation}\label{SpNorm-In}
\|A{\bf x}\|_2\leq \|A\|_2\|{\bf x}\|_2
\end{equation}
for any appreciable ${\bf x}\in \widehat{\mathbb Q}^n$.
The following proposition shows that, similar to the common complex matrix situation,  $\|A\|_2$ defined by \eqref{SNorm-DQM} is exactly the largest singular value of $A$. 

 \begin{Prop}
 Let $A\in \widehat{\mathbb Q}^{m\times n}$. It holds that $\|A\|_2=\sigma_1(A)$, where $\sigma_1(A)$ is the largest singular value of $A$.
 \end{Prop}

 \begin{proof}
Let $A=U\Sigma V^*$ be a singular value decomposition of $A$, in which $U$ and $V$ are unitary, $\Sigma={\rm diag}(\sigma_1(A),\ldots,\sigma_s(A))$ with $\sigma_1(A)\geq\ldots\geq\sigma_s(A)\geq0$ and $s=\min\{m,n\}$. It follows from (\ref{SNorm-DQM}) and Proposition \ref{pp3.2} that
$$
\begin{array}{lll}
\|A \|_2 &=&\displaystyle \max_{{\bf x}\in \widehat{\mathbb Q}^n,~\|{\bf x}\|_2=1}\|\Sigma V^*{\bf x}\|_2\\
&=&\displaystyle \max_{{\bf y}\in \widehat{\mathbb Q}^n,~\|V{\bf y}\|_2=1}\|\Sigma {\bf y}\|_2\\
&=&\displaystyle \max_{{\bf y}\in \widehat{\mathbb Q}^n,~\|{\bf y}\|_2=1}\sqrt{\sum_{i=1}^n\sigma^2_i(A) |y_i|^2}\\
&\leq&\displaystyle \max_{{\bf y}\in \widehat{\mathbb Q}^n,~\|{\bf y}\|_2=1}\sigma_1(A) \sqrt{\sum_{i=1}^n|y_i|^2}\\
&=&\sigma_1(A),
\end{array}
$$
where the inequality comes from items (a), (c) and (f) of Proposition \ref{P6.5}. However, $\|A{\bf y}\|_2= \sigma_1(A)$ for ${\bf y}={\bf e}_1$, hence $\|A \|_2\geq \sigma_1(A)$. Therefore, we conclude that $\|A \|_2= \sigma_1(A)$ and complete the proof.
\end{proof}

\section{von Neumann type trace inequality}\label{VN-TraceIn}
In this section, we extend the well-known von Neumann trace inequality to dual quaternionic versions. Because $\widehat{\mathbb R}$ is a total order space in the meaning of total order stated in Section \ref{Du-number}, unless otherwise specified, we have, for $p,q\in \widehat{\mathbb R}$, $p\leq q$ if and only if $q-p\in \widehat{\mathbb R}_+$. We start this section by introducing the following concept for dual numbers.

For given dual numbers $x_1,x_2,\ldots, x_m$, we use $\hat x_1,\hat x_2,\ldots, \hat x_m$ to denote these numbers arranged in non-ascending order of magnitude. If the two sets of dual numbers $x_1,x_2,\ldots, x_m$ and $y_1,y_2,\ldots, y_m$ satisfy the
relations
$$
\hat x_1+\hat x_2+\ldots+ \hat x_s\left\{\begin{array}{ll}
\leq\hat y_1+\hat y_2+\ldots+ \hat y_s&{\rm for~}1\leq s\leq m-1,\\
=\hat y_1+\hat y_2+\ldots+ \hat y_s&{\rm for ~}s=m,
\end{array}\right.
$$
we write $(x_1,x_2, \ldots,x_m)\prec (y_1, y_2, \ldots, y_m)$ for simplicity.

The following lemma is an extension of Lemma in \cite{Mir59} to the case of dual numbers. Due to the introduction of the total order ``$\geq$'' in $\widehat{\mathbb R}$, it has similar properties to the set of real numbers. For example, if $a,b\in \widehat{\mathbb R}_+$ implies $a+b\in \widehat{\mathbb R}_+$ and $ab\in \widehat{\mathbb R}_+$. Although the proof of this lemma is similar to the way used for the unique lemma in \cite{Mir59}, we give a proof for the completeness of this paper.
\begin{Lem}\label{xyzIn}
Let $\{x_1,x_2,\ldots, x_m\}$, $\{y_1,y_2,\ldots, y_m\}$ and $\{z_1,z_2,\ldots, z_m\}\subset \widehat{\mathbb R}$. Suppose $x_1\geq x_2\geq \ldots\geq x_m$, $y_1\geq y_2\geq \ldots\geq y_m$ and $(z_1,z_2, \ldots,z_m)\prec (y_1, y_2, \ldots, y_m)$. Then it holds that
\begin{equation}
\sum_{i=1}^mx_iz_i\leq \sum_{i=1}^mx_iy_i.
\end{equation}
\end{Lem}
\begin{proof}
For any $1\leq k\leq m$, let $X_k=x_1+ x_2+\ldots+ x_k$ and $Z_k=\hat z_1+\hat z_2+\ldots+ \hat z_k$. Then, by virtue of hypothesis, we have $Z_k\leq Y_k$ for $k=1,2,\ldots,m$. Consequently, it holds that
$$
\begin{array}{lll}
\displaystyle\sum_{k=1}^mx_kz_k&\leq&\displaystyle\sum_{k=1}^mx_k\hat z_k\\
&=&\displaystyle x_1Z_1+\sum_{k=2}^mx_k(Z_k-Z_{k-1})\\
&=&\displaystyle\sum_{k=1}^{m-1}(x_k-x_{k+1})Z_k+x_nZ_n\\
&\leq&\displaystyle\sum_{k=1}^{m-1}(x_k-x_{k+1})Y_k+x_nY_n\\
&=&\displaystyle\sum_{k=1}^mx_ky_k.
\end{array}
$$
The proof is completed.
\end{proof}

\begin{Lem}\label{Lemma2}
Let $A=(a_{ij})\in \widehat{\mathbb Q}^{m\times m}$ be a Hermitian matrix, and let the right eigenvalues $\lambda_i(A)$ of $A$ be arranged so that $\lambda_1(A)\geq \lambda_2(A)\geq \ldots\geq \lambda_m(A)$. Then, for any given positive integer $k\leq m$, the sum $\sum_{i=1}^k\lambda_i(A)$ is the maximum of $\sum_{i=1}^k({\bf x}^{(i)})^*A{\bf x}^{(i)}$, when $k$ orthonormal vectors ${\bf x}^{(i)}~(i=1,2,\ldots,k)$ vary in $\widehat{\mathbb Q}^m$. In particular, we have
\begin{equation}\label{TrEigIN}
\sum_{i=1}^k a_{ii}\leq \sum_{i=1}^k\lambda_i(A), ~~~~k=1,2,\ldots,m.
\end{equation}
\end{Lem}

\begin{proof}
By Theorem \ref{HUDec}, there exists a unitary matrix $U\in \widehat{\mathbb Q}^{m\times m}$ such that
$$U^*AU={\rm diag}(\lambda_1(A),\lambda_2(A),\ldots,\lambda_m(A)),$$
which implies $A{\bf u}^{(i)}={\bf u}^{(i)}\lambda_i(A)$ for $i=1,2,\ldots,m$, where ${\bf u}^{(i)}$ is the $i$th column of $U$. Since $UU^*=I$,  it holds, for any ${\bf x}^{(j)}\in \widehat{\mathbb Q}^m$, that
$$
{\bf x}^{(j)}=UU^*{\bf x}^{(j)}=\sum_{i=1}^m{\bf u}^{(i)}(({\bf u}^{(i)})^*{\bf x}^{(j)})=\sum_{i=1}^m{\bf u}^{(i)}\langle{\bf x}^{(j)},{\bf u}^{(i)}\rangle.
$$
Since $A\in \widehat{\mathbb Q}^{m\times m}$ is a Hermitian matrix, all $\lambda_i$'s are dual numbers. Consequently, we have
$$
\begin{array}{lll}
({\bf x}^{(j)})^*A{\bf x}^{(j)}&=&\displaystyle\sum_{i=1}^m\lambda_i(A)|\langle{\bf x}^{(j)},{\bf u}^{(i)}\rangle|^2\\
&=&\displaystyle\lambda_k(A)\sum_{i=1}^m|\langle{\bf x}^{(j)},{\bf u}^{(i)}\rangle|^2+\sum_{i=1}^k(\lambda_i(A)-\lambda_k(A))|\langle{\bf x}^{(j)},{\bf u}^{(i)}\rangle|^2\\
&&+\displaystyle\sum_{i=k+1}^m(\lambda_i(A)-\lambda_k(A))|\langle{\bf x}^{(j)},{\bf u}^{(i)}\rangle|^2.
\end{array}
$$
Hence, when $\|{\bf x}^{(j)}\|_2=1$, we obtain
$$
({\bf x}^{(j)})^*A{\bf x}^{(j)}\leq \lambda_k(A)+\sum_{i=1}^k(\lambda_i(A)-\lambda_k(A))|\langle{\bf x}^{(j)},{\bf u}^{(i)}\rangle|^2,
$$
since
$$
\sum_{i=1}^m|\langle{\bf x}^{(j)},{\bf u}^{(i)}\rangle|^2=({\bf x}^{(j)})^*UU^*{\bf x}^{(j)}=({\bf x}^{(j)})^*{\bf x}^{(j)}=\|{\bf x}^{(j)}\|_2=1,
$$
as well as $\lambda_i(A)\leq \lambda_k(A)$ for $i=k+1,\ldots, m$. By this, it holds that
$$
\begin{array}{lll}
\displaystyle\sum_{j=1}^k({\bf x}^{(j)})^*A{\bf x}^{(j)}&\leq &\displaystyle\sum_{j=1}^k\lambda_k(A)+\sum_{j=1}^k\sum_{i=1}^k(\lambda_i(A)-\lambda_k(A))|\langle{\bf x}^{(j)},{\bf u}^{(i)}\rangle|^2\\
&=&\displaystyle\sum_{j=1}^k\lambda_k(A)+\sum_{i=1}^k(\lambda_i(A)-\lambda_k(A))\sum_{j=1}^k|\langle{\bf x}^{(j)},{\bf u}^{(i)}\rangle|^2,
\end{array}
$$
which implies
$$
\begin{array}{lll}
\displaystyle\sum_{i=1}^k\lambda_i(A)-\sum_{j=1}^k({\bf x}^{(j)})^*A{\bf x}^{(j)}&\geq &\displaystyle\sum_{j=1}^k(\lambda_j(A)-\lambda_k(A))-\sum_{i=1}^k\sum_{j=1}^k\left(\lambda_i(A)-\lambda_k(A)\right)|\langle{\bf x}^{(j)},{\bf u}^{(i)}\rangle|^2\\
&=&\displaystyle\sum_{j=1}^k(\lambda_j(A)-\lambda_k(A))\left\{1-\sum_{i=1}^k|\langle{\bf x}^{(j)},{\bf u}^{(i)}\rangle|^2\right\}\\
&\geq &0,
\end{array}
$$
since $\lambda_j(A)\geq\lambda_k(A)$ for $j=1,2,\ldots,k$ and
$$
\sum_{i=1}^k|\langle{\bf x}^{(j)},{\bf u}^{(i)}\rangle|^2=({\bf x}^{(j)})^*\widetilde{U}\widetilde{U}^*{\bf x}^{(j)}\leq \|\widetilde{U}\widetilde{U}^*\|_2\|{\bf x}^{(j)}\|^2_2\leq\|{\bf x}^{(j)}\|^2_2=1,
$$
where $\widetilde{U}=[{\bf u}^{(1)},\ldots,{\bf u}^{(k)}]$. On the other hand, by taking ${\bf x}^{(j)}={\bf u}^{(j)}$ for $j=1,2,\ldots,k$, it is easy to verify that  $\sum_{i=1}^k\lambda_i(A)=\sum_{j=1}^k({\bf x}^{(j)})^*A{\bf x}^{(j)}$.

Notice that $a_{ii}\in \widehat{\mathbb R}$ for $i=1,2,\ldots,m$, since $A$ is Hermitian. In particular, by taking ${\bf x}^{(j)}={\bf e}_j$, where ${\bf e}_j$ is the $j$th column in $I_m$, it is easy to see that (\ref{TrEigIN}) holds. Therefore, we obtain the desired result and complete the proof.
\end{proof}

Hereafter, we recall a lemma that can be found in \cite[Theorem 4.10]{LQY22}.
\begin{Lem}[\cite{LQY22}]\label{Lemma3}
Let $A\in \widehat{\mathbb Q}^{m\times m}$ and $H=(A^*+A)/2$. Let $\sigma_1(A)\geq\sigma_2(A)\geq\ldots\geq\sigma_m(A)$ be singular values of $A$, and let $\lambda_1(H)\geq\lambda_2(H)\geq\ldots\geq\lambda_m(H)$ be the right eigenvalues  of $H$. Then, it holds that
$$
\lambda_i(H)\leq\sigma_i(A), ~~~i=1,2,\ldots,m.
$$
\end{Lem}

With the above preparations, we now state and prove the von Neumann type trace inequality for dual quaternion matrices.  Here, we prove the dual quaternionic von Neumann inequality by a unified way, which is applicable to the quaternionic case and is completely different with the proof given in \cite{CY16}.

\begin{Thm}\label{VNINTH}
For any $A,B\in \widehat{\mathbb Q}^{m\times n}$, it holds that
\begin{equation}
{\rm trace}(A^*B+B^*A)\leq 2\sum_{i=1}^s\sigma_i(A)\sigma_i(B),
\end{equation}
where $s=\min\{m,n\}$, and $\sigma_1(A)\geq\sigma_2(A)\geq\ldots\geq\sigma_s(A)$ and $\sigma_1(B)\geq\sigma_2(B)\geq\ldots\geq\sigma_s(B)$ are singular values of $A$ and $B$, respectively.
\end{Thm}

\begin{proof}
Without loss of generality, we assume $m\leq n$. Let $A=U\Sigma_AV^*$ be the SVD of $A$, where $\Sigma_A={\rm diag}(\sigma_1(A),\sigma_2(A),\ldots,\sigma_m(A))$. It is obvious that $A^*=V\Sigma_A^\top U^*$, since all $\sigma_i(A)$ ($i=1,2,\ldots,m$) are dual numbers. Consequently, we have
$$
V^*(A^*B+B^*A)V=\Sigma_A^\top C+C^* \Sigma_A,
$$
where $C=U^*BV\in \widehat{\mathbb Q}^{m\times n}$. It is easy to see that
$$
{\rm trace}(A^*B+B^*A)={\rm trace}\left\{V^*(A^*B+B^*A)V\right\}={\rm trace}(\Sigma_A^\top C+C^* \Sigma_A)=\sum_{i=1}^m\sigma_i(A)(c_{ii}+\bar c_{ii}),
$$
where $c_{ii}$ is the $i$th diagonal element of $C$ and $\bar c_{ii}$ is the conjugate of the dual quaternion number $c_{ii}$.

Now, we prove
$$
\sum_{i=1}^m\sigma_i(A)(c_{ii}+\bar c_{ii})\leq2\sum_{i=1}^{m}\sigma_i(A)\sigma_i(B).
$$
Let $$D=(d_{ij}):=\left[\begin{array}{c}C\\
0_{(n-m)\times n}\end{array}\right]+[C^*,0_{n\times (n-m)}].$$ It is obvious that $D\in \widehat{\mathbb Q}^{n\times n}$ satisfies $D^*=D$, i.e., $D$ is a dual quaternion Hermitian matrix. It is obvious that $d_{ii}=c_{ii}+\bar c_{ii}$ for any $i=1,2,\ldots,m$. Applying Lemma \ref{Lemma2} with $A=D$, we obtain 
\begin{equation}\label{ccdSinIn}
\sum_{i=1}^k(c_{ii}+\bar c_{ii})=\sum_{i=1}^kd_{ii}\leq \sum_{i=1}^k\lambda_i(D)\leq2\sum_{i=1}^k\sigma_i(C)=2\sum_{i=1}^k\sigma_i(B), ~~k=1,2,\ldots,m,
\end{equation}
where the second inequality comes from Lemma \ref{Lemma3} and the fact that $\sigma_i\left(\left[\begin{array}{c}C\\0_{(n-m)\times m}\end{array}\right]\right)=\sigma_i(C)$ for $i=1,2,\ldots,m$, and the last equality is due to the fact $\sigma_i(C)=\sigma_i(B)$ for $i=1,2,\ldots,m$ since $U,V$ are unitary. Finally, since $\sigma_i(A)\geq $ for $i=1,2,\ldots, m$, by (\ref{ccdSinIn}) and  Lemma \ref{xyzIn}, we obtain the desired conclusion and complete the proof.
\end{proof}

\begin{Prop}
Let $A\in \widehat{\mathbb Q}^{m\times n}$. We have
$$
\mathop{\rm max}\limits_{X\in \mathbb{U}\widehat{\mathbb Q}^{m\times n}}{\rm trace}(A^*X+X^*A)= 2\sum_{i=1}^s\sigma_i(A),
$$
where $\mathbb{U}\widehat{\mathbb Q}^{m\times n}=\{X\in \widehat{\mathbb Q}^{m\times n}~|~XX^*=I_m\}$ and $\sigma_i(A)$ is the $i$-th singular value of $A$ for $1\leq i\leq s={\min}\{m,n\}$.
\end{Prop}

\begin{proof}
Without loss of generality, we assume $m\leq n$. Let $A=U\Sigma_AV^*$ be the SVD of $A$, where $\Sigma_A={\rm diag}(\sigma_1(A),\sigma_2(A),\ldots,\sigma_m(A))$. Denote $\bar X=USV^*$, where $S=[E_m,O_{m\times (n-m)}]\in \widehat{\mathbb Q}^{m\times n}$. It is obvious that $\bar X\bar X^*=I_m$ and $\sigma_i(\bar X)=1$ for $i=1,2,\ldots, m$. Moreover, it is clear that $A^*=V\Sigma_A^\top U^*$, since all $\sigma_i(A)$ ($i=1,2,\ldots,m$) are dual numbers. Consequently, it is easy to see that
$$
{\rm trace}(A^*\bar X+\bar X^*A)=2{\rm trace}\left(V\left[\begin{array}{cc}
\Sigma_A&O_{m\times (n-m)}\\
O_{(n-m)\times m}&O_{(n-m)\times (n-m)}
\end{array}\right]V^*\right)=2{\rm trace}(\Sigma_A)=2\sum_{i=1}^m\sigma_i(A),
$$
where the second equality follows from the fact that $\left[\begin{array}{cc}
\Sigma_A&O_{m\times (n-m)}\\
O_{(n-m)\times m}&O_{(n-m)\times (n-m)}
\end{array}\right]$ is Hermitian. Finally, by combining the above discussion and Theorem \ref{VNINTH}, we obtain the desired result and complete the proof.
\end{proof}

Due to the application of dual quaternion Hermitian matrices in multi-agent formation control, we below present a variant of the von Neumann inequality expressed by eigenvalues of dual quaternion Hermitian matrices.

\begin{Thm}\label{EigFnormIn}
	Let $A,B \in \widehat{\mathbb{Q}}^{m\times m}$. Suppose that both $A$ and $B$ are Hermitian. We have
	$$
	{\rm trace}(AB+BA)\leq 2\sum_{i=1}^m\lambda_i(A)\lambda_i(B),
	$$
	where $\lambda_1(A)\geq\lambda_2(A)\geq\ldots\geq\lambda_m(A)$ and $\lambda_1(B)\geq\lambda_2(B)\geq\ldots\geq\lambda_{m}(B)$ are the eigenvalues of $A$ and $B$, respectively. In particular, the above inequality holds when $A,B \in \mathbb{Q}^{m\times m}$ are Hermitian.
\end{Thm}

\begin{proof}
	Since both $A$ and $B$ are dual quaternion Hermitian matrices, we know that $\lambda_i(A), \lambda_i(B)\in \widehat{\mathbb{R}}$ for $i=1,2,\ldots,m$. By Theorem \ref{HUDec}, there exists a unitary matrix $U\in \widehat{\mathbb Q}^{m\times m}$ such that $A=U\Sigma U^*$ where $\Sigma={\rm diag} (\lambda_1(A),\lambda_2(A),\ldots,\lambda_m(A))$. Let $C:=(c_{ij})=U^*BU$. It is clear that $C$ is a dual quaternion Hermitian matrix, which implies $c_{ii}\in \widehat{\mathbb{R}}$ for $i=1,2,\ldots,m$. Moreover, we have
	\begin{equation}\label{trABEq}
	{\rm trace}(AB+BA)={\rm trace}(U(\Sigma C+C\Sigma)U^*)={\rm trace}(\Sigma C+C\Sigma)=2\sum_{i=1}^m\lambda_i(A)c_{ii}.
	\end{equation}
	Since $C$ is Hermitian, by Lemma \ref{Lemma2} and {\color {blue} (\ref{sumEigV})}, for every $k=1,2,\ldots m-1$, it holds that
	$$
	\sum_{i=1}^kc_{ii}\leq \sum_{i=1}^k\lambda_i(C)=\sum_{i=1}^k\lambda_i(B)~~{\rm and}~~\sum_{i=1}^mc_{ii}=\sum_{i=1}^m\lambda_i(B),
	$$
	since $\lambda_i(C)=\lambda_i(B)$ for $i=1,2,\ldots,m$. Consequently, by Lemma \ref{xyzIn}, we have  $$
	\sum_{i=1}^m\lambda_i(A)c_{ii}\leq \sum_{i=1}^m\lambda_i(A)\lambda_i(B),
	$$
	which implies, together with (\ref{trABEq}), the desired inequality holds.
\end{proof}

\section{Hoffman-Wielandt type inequality}\label{sec_HWineq}

In this section, we are concerned with the Hoffman-Wielandt type inequality for dual quaternion matrices.

First, by Theorem \ref{VNINTH}, we have the following theorem, which is also a generalization of the well-known Hoffman-Wielandt type inequality \cite{HW53} in $\mathbb{C}^{m\times n}$, and characterizes an upper bound for all singular values simultaneous perturbation of a dual quaternion matrix.

\begin{Thm}\label{FNormIn}
	Let $A,B\in \widehat{\mathbb Q}^{m\times n}$. If $A-B$ is appreciable, 
	then it holds that
	\begin{equation}\label{FNormIn-1}
	\|\sigma(A)-\sigma(B)\|_2\leq \|A-B\|_F,
	\end{equation}
	where $s=\min\{m,n\}$, $\sigma(A)=(\sigma_1(A), \ldots,\sigma_{s}(A))^\top$ and $\sigma(B)=(\sigma_1(B), \ldots,\sigma_{s}(B))^\top$ with $\sigma_1(A)\geq\sigma_2(A)\geq\ldots\geq\sigma_{s}(A)$ and $\sigma_1(B)\geq\sigma_2(B)\geq\ldots\geq\sigma_{s}(B)$ being the singular values of $A$ and $B$, respectively.
\end{Thm}

\begin{proof}
	Let $A=U\Sigma_AV^*$ and $B=X\Sigma_BY^*$ be the SVDs of $A$ and $B$, respectively. Under the condition $A-B$ being appreciable, we divide our proofs into two cases.
	\begin{itemize}
		\item If $\sigma(A)-\sigma(B)$ is infinitesimal, then $A-B$ being appreciable implies
		$$\|A-B\|_F=\|A_{\sf st}-B_{\sf st}\|_F+\delta\epsilon$$
		for some $\delta\in \widehat{\mathbb R}_+$ by Proposition \ref{p6.3}. Clearly, it follows from the fact $\|A_{\sf st}-B_{\sf st}\|_F>0$ that (\ref{FNormIn-1}) holds.
		\item If $\sigma(A)-\sigma(B)$ is appreciable, then we only need to prove
		$$
		\sum_{i=1}^{s}|\sigma_i(A)-\sigma_i(B)|^2\leq \|A-B\|^2_F,
		$$
		which is equivalent to
		\begin{equation}\label{FNormIn-2}
		\sum_{i=1}^{s}\sigma^2_i(A)-2\sum_{i=1}^{s}\sigma_i(A)\sigma_i(B)+\sum_{i=1}^{s}\sigma_i(B)^2\leq \|A\|_F^2-{\rm trace}(A^*B+B^*A)+\|B\|^2_F.
		\end{equation}
		Recalling the fact that $\sum_{i=1}^{s}\sigma^2_i(A)=\|A\|_F^2$ and $\sum_{i=1}^{s}\sigma_i(B)^2=\|B\|^2_F$,  we immediately prove (\ref{FNormIn-2}) with the employment of Theorem \ref{VNINTH}.
	\end{itemize}
	To sum up, we obtain the desired conclusion and complete the proof.
\end{proof}

\begin{ReK}
	Notice that, if $A,B\in \mathbb{Q}^{m\times n}$, it is obvious that Theorem \ref{VNINTH} holds. Therefore, if both $A$ and $B$  are either quaternions or infinitesimal, Theorem \ref{FNormIn} holds. In fact, in the case that both $A$ and $B$  are infinitesimal, we claim that $\sigma(A)$ and $\sigma(B)$ are both infinitesimal, and $(\sigma(A))_{\sf in}=\sigma(A_{\sf in})$ and $(\sigma(B))_{\sf in}=\sigma(B_{\sf in})$. Consequently, we have
	$$\|\sigma(A)-\sigma(B)\|_2=\|(\sigma(A))_{\sf in}-(\sigma(B))_{\sf in}\|_2\epsilon=\|\sigma(A_{\sf in})-\sigma(B_{\sf in})\|_2\epsilon$$ and $$\|A-B\|_F=\|A_{\sf in}-B_{\sf in}\|_F\epsilon.$$
	Hence, we only need to prove $\|\sigma(A_{\sf in})-\sigma(B_{\sf in})\|_2\leq \|A_{\sf in}-B_{\sf in}\|_F$, which can be proved by applying Theorem \ref{VNINTH} with $A=A_{\sf in}$ and $B=B_{\sf in}$. The conclusion for the case $A, B\in \mathbb{Q}^{m\times n}$ can be proved similarly.
\end{ReK}

\begin{ReK}\label{remark3}
	From Theorem \ref{SVD-DQM}, we know that the standard parts of the singular values of a dual quaternion matrix are exactly the singular values of the standard part of that dual quaternion matrix. Hence, If both $A$ and $B$ are appreciable, but $A-B$ is infinitesimal, i.e., $A_{\sf st}=B_{\sf st}\neq O$, then $\sigma(A)-\sigma(B)$ is infinitesimal. In this case, the desired inequality \eqref{FNormIn-1} becomes $\|(\sigma(A)-\sigma(B))_{\sf in}\|_2\leq \|A_{\sf in}-B_{\sf in}\|_F$. However, we do not know whether this inequality still holds, and leave it as an open question for one of our future concerns.
\end{ReK}

As mentioned in Remark \ref{remark3}, the infinitesimal part of dual numbers makes the analysis on dual quaternion matrices is difficult, and Theorem \ref{FNormIn} holds under the condition that $A-B$ is appreciable. Below, we are concerned with whether the Hoffman-Wielandt inequality still holds for dual quaternion Hermitian matrices when removing the condition $A-B$ being appreciable.

We first show the quaternionic Hoffman-Wielandt inequality.

\begin{Prop}\label{SeProp2}
Let $A,B \in \mathbb{Q}^{m\times m}$. Suppose that both $A$ and $B$ are Hermitian. We have
$$
\|\lambda(A)-\lambda(B)\|_2\leq \|A-B\|_F,
$$
where $\lambda(A)=(\lambda_1(A), \ldots,\lambda_m(A))^\top$ and $\lambda(B)=(\lambda_1(B), \ldots,\lambda_m(B))^\top$ with $\lambda_1(A)\geq\lambda_2(A)\geq\ldots\geq\lambda_m(A)$ and $\lambda_1(B)\geq\lambda_2(B)\geq\ldots\geq\lambda_{m}(B)$ being the eigenvalues of $A$ and $B$, respectively.
\end{Prop}

\begin{proof}
Since ${\rm trace}(A^2)=\|A\|_F^2=\|\lambda(A)\|_2^2$ and ${\rm trace}(B^2)=\|B\|_F^2=\|\lambda(B)\|_2^2$, it follows from Proposition \ref{EigFnormIn}.
\end{proof}

Now we state and prove the following theorem, which indicates that the Hoffman-Wielandt type inequality for dual quaternion Hermitian matrices still holds even if $A-B$ is  infinitesimal.

\begin{Thm}
Let $A, B\in \widehat{\mathbb{Q}}^{m\times m}$. If both $A$ and $B$ are Hermitian matrices, then we have
\begin{equation}\label{HW-Ineq}
\|\lambda(A)-\lambda(B)\|_2\leq \|A-B\|_F,
\end{equation}
where $\lambda(A)=(\lambda_1(A), \ldots,\lambda_m(A))^\top$ and $\lambda(B)=(\lambda_1(B), \ldots,\lambda_m(B))^\top$ with $\lambda_1(A)\geq\lambda_2(A)\geq\ldots\geq\lambda_m(A)$ and $\lambda_1(B)\geq\lambda_2(B)\geq\ldots\geq\lambda_{m}(B)$ being the eigenvalues of $A$ and $B$, respectively.
\end{Thm}

\begin{proof}
Firstly, since both $A$ and $B$ are Hermitian, we know that $\lambda(A), \lambda(B)\in \widehat{\mathbb{R}}^m$. We divide the proofs into two cases: (a) $A-B$ is appreciable, and (b) $A-B$ is infinitesimal.

In the case (a), since $\|\lambda(A)\|_2^2=\|A\|_F^2$ and $\|\lambda(B)\|_2^2=\|B\|_F^2$, regardless of whether $\lambda(A)-\lambda(B)$ is appreciable or not, by Proposition \ref{EigFnormIn}, we can prove (\ref{HW-Ineq})  in a similar way to the proof of Theorem \ref{FNormIn}.

We now prove the desired inequality (\ref{HW-Ineq}) for the case where $A-B$ is infinitesimal. In this case, both $A$ and $B$ must be appreciable or infinitesimal at the same time, since $A_{\sf st}=B_{\sf st}$. If both $A$ and $B$ are infinitesimal, i.e., $A_{\sf st}=B_{\sf st}=O$, then the inequality (\ref{HW-Ineq}) becomes
\begin{equation}\label{HW-InQ}
\|(\lambda(A))_{\sf in}-(\lambda(B))_{\sf in}\|_2\leq \|A_{\sf in}-B_{\sf in}\|_F,
\end{equation} since $\lambda(A)=\lambda(A_{\sf in})\epsilon$ and $\lambda(B)=\lambda(B_{\sf in})\epsilon$. Consequently, since both $A_{\sf in}$ and $B_{\sf in}$ are Hermitian, by Proposition \ref{EigFnormIn} again, we know that the inequality (\ref{HW-InQ}) holds.

If $A$ and $B$ are appreciable Hermitian matrices at the same time, i.e., $A_{\sf st}=B_{\sf st}$ is a nonzero quaternion Hermitian matrix, then by quaternion matrix theory, there exists a unitary matrix $S\in \mathbb{Q}^{m\times m}$, such that
$$
SAS^*=D+C\epsilon~~~~{\rm and}~~~~SBS^*=D+G\epsilon,
$$
where $C=SA_{\sf in}S^*$, $G=SB_{\sf in}S^*$, $D =SA_{\sf st}S^*={\rm diag}(\lambda_1I_{k_1}, \lambda_2I_{k_2}, \ldots, \lambda_rI_{k_r})$. Here, $\lambda_1 > \lambda_2 > \cdots > \lambda_r$ are real numbers, $I_{k_i}$ is a $k_i \times k_i$ identity matrix, and $\sum_{i=1}^r k_i = m$. It is obvious that
$\|A_{\sf in}-B_{\sf in}\|_F^2=\|C-G\|_F^2$. Notice that both $C$
and $G$ are quaternion Hermitian matrices. Let
\begin{equation*}
C=\begin{bmatrix}
C_{11} &  C_{12} & \cdots &   C_{1r} \\
  C_{12}^* & C_{22} & \cdots  &   C_{2r} \\
\vdots & \vdots & \ddots & \vdots \\
  C_{1r}^* &   C_{2r}^* &  \cdots &  C_{rr}
\end{bmatrix}
~~~~{\rm and}~~~~
G=\begin{bmatrix}
G_{11} &  G_{12} & \cdots &   G_{1r}\\
  G_{12}^*& G_{22} & \cdots  &   G_{2r} \\
\vdots & \vdots & \ddots & \vdots \\
  G_{1r}^*&   G_{2r}^* &  \cdots &  G_{rr}
\end{bmatrix},
\end{equation*}
where $C_{ij}$ and $G_{ij}$ are quaternion matrices of same adequate dimensions, and $C_{ii}$ and $G_{ii}$ are Hermitian for $i=1,2,\ldots,r$. Since $C_{ii}$ and $G_{ii}$ are Hermitian, there exist real numbers $\lambda_{i,1} \ge \ldots \ge \lambda_{i,k_i}$ and $\mu_{i,1} \ge \ldots \ge \mu_{i, k_i}$ for $i=1,2,\ldots,r$,  such that
\begin{equation}\label{blocks}
C_{ii} = U_{i}{\rm diag}\left(\lambda_{i,1}, \cdots, \lambda_{i,k_i}\right)U_{i}^* ~~{\rm and}~~G_{ii}= V_{i} {\rm diag}\left(\mu_{i,1}, \cdots, \mu_{i,k_i}\right)V_{i}^* , ~~i=1,\cdots, r.
\end{equation}
It is obvious that $\lambda_{i,j}$ and $\mu_{i,j}$ are the $j$th largest eigenvalues of $C_{ii}$ and $G_{ii}$ respectively, for $i=1,2,\ldots,r$. From the proof precess of Theorem 4.1 presented in \cite{QL21}, we know that
$$
\lambda(A)_{\sf in}=(\lambda_{1,1},\ldots,\lambda_{1,k_1},\lambda_{2,1},\ldots,\lambda_{2,k_2},\ldots,\lambda_{r,1},\ldots,\lambda_{r,k_r})^\top
$$
and
$$
\lambda(B)_{\sf in}=(\mu_{1,1},\ldots,\mu_{1,k_1},\mu_{2,1},\ldots,\mu_{2,k_2},\ldots,\mu_{r,1},\ldots,\mu_{r,k_r})^\top.
$$
Consequently, it holds that
$$
\begin{array}{lll}
\|\lambda(A)_{\sf in}-\lambda(B)_{\sf in}\|_2^2&=&\displaystyle\sum_{i=1}^r\sum_{j=1}^{k_i}|\lambda_{i,j}-\mu_{i,j}|^2\\
&\leq &\displaystyle\sum_{i=1}^r\|C_{ii}-G_{ii}\|^2\\
&\leq&\displaystyle\sum_{i=1}^r\sum_{j=1}^r\|C_{ij}-G_{ij}\|^2\\
&=&\|C-G\|_F^2\\
&=&\|A_{\sf in}-B_{\sf in}\|_F^2,
\end{array}
$$
where the first inequality comes from Proposition \ref{SeProp2} with $A=C_{ii}$ and $B=G_{ii}$ for $i=1,2,\ldots,r$. Hence, we obtain $\|\lambda(A)_{\sf in}-\lambda(B)_{\sf in}\|_2\leq\|A_{\sf in}-B_{\sf in}\|_F$, which implies $\|\lambda(A)-\lambda(B)\|_2\leq\|A-B\|_F$ by (\ref{FNorm-DQM}), since $\lambda(A)_{\sf st}=\lambda(B)_{\sf st}$ and $A_{\sf st}=B_{\sf st}$. We obtain the desired result and complete the proof.
\end{proof}

\section{Conclusion}\label{Conclusion}
In this paper, after introducing the concept of spectral norm for dual quaternion matrices, we extended the well-known von Neumann trace inequality for general dual quaternion matrices. Using the proposed trace inequality, we further obtained a Hoffman-Wielandt type inequality, which characterizes the distance between two dual quaternion matrices $A, B\in \widehat{\mathbb Q}^{m\times n}$ being larger than the distance between their respective singular values. Such an inequality can also be regarded as a simultaneous perturbation bound on all singular values of a general dual quaternion matrix. Furthermore, we also proposed two variants of the above two inequalities expressed by eigenvalues of dual quaternion Hermitian matrices. In particular, we proved that the Hoffman-Wielandt type inequality still holds even if $A-B$ is infinitesimal. As a new area of applied mathematics, there are many problems worth exploring to enrich the theory of dual quaternion matrices, such as optimal low-rank approximations and applications of dual quaternion matrices in the fields of data analysis, computer science and intelligent control.

\begin{thebibliography}{10}
	
	\bibitem{BK20}
	G.~Brambley and J.~Kim.
	\newblock Unit dual quaternion-based pose optimization for visual runway
	observations.
	\newblock {\em IET Cyber Systems and Robotics}, 2:181--189, 2020.
	
	\bibitem{BLH19}
	S.~Bultmann, K.~Li, and U.D. Hanebeck.
	\newblock Stereo visual {SLAM} based on unscented dual quaternion filtering.
	\newblock In {\em 2019 22th International Conference on Information Fusion
		(FUSION)}, pages 1--8, 2019.
	
	\bibitem{CY16}
	T.-S. Chan and Y.-H. Yang.
	\newblock Complex and quaternion principal component pursuit and its
	application to audio separation.
	\newblock {\em IEEE Signal Processing Letters}, 23:287--291, 2016.
	
	\bibitem{CQZ22}
	Y.~Chen, L.~Qi, and X.~Zhang.
	\newblock Color image completion using a low-rank quaternion matrix
	approximation.
	\newblock {\em Pacific Journal of Optimization}, 18:55--75, 2022.
	
	\bibitem{CXZ20}
	Y.~Chen, X.~Xiao, and Y.~Zhou.
	\newblock Low-rank quaternion approximation for color image processing.
	\newblock {\em IEEE Transactions on Image Processing}, 29:1426--1439, 2020.
	
	\bibitem{CKJC16}
	J.~Cheng, J.~Kim, Z.~Jiang, and W.~Che.
	\newblock Dual quaternion-based graph {SLAM}.
	\newblock {\em Robotics and Autonomous Systems}, 77:15--24, 2016.
	
	\bibitem{Cl73}
	W.~K. Clifford.
	\newblock Preliminary sketch of bi-quaternions.
	\newblock {\em Proceedings of the London Mathematical Society}, 4:381--395,
	1873.

    \bibitem{CQ23}
	C. Cui and L. Qi.
	\newblock A power method for computing the dominant eigenvalue of a dual quaternion Hermitian matrix.
    \newblock arXiv: 2304.04355, 2023.
	
	
	\bibitem{Da99}
	K.~Daniilidis.
	\newblock Hand-eye calibration using dual quaternions.
	\newblock {\em The International Journal of Robotics Research}, 18:286--298,
	1999.
	
	\bibitem{Gri17}
	S.~Griffin.
	\newblock {\em Quaternions: Theory and Applications}.
	\newblock Mathematics Research Developments. NOVA Science Publishers, New York,
	2017.
	
	\bibitem{HW53}
	A.~J. Hoffman and H.~W. Wielandt.
	\newblock The variation of the spectrum of a normal matrix.
	\newblock {\em Duke Mathematical Journal}, 20:37--40, 1953.
	
	\bibitem{HJ13}
	R.A. Horn and C.R. Johnson.
	\newblock {\em Matrix Analysis}.
	\newblock Cambridge University Press, New York, second edition edition, 2013.
	
	\bibitem{Ke12}
	B.~Kenwright.
	\newblock A beginners guide to dual-quaternions.
	\newblock In {\em The 20th International Conference on Computer Graphics,
		Visualization and Computer Vision}, pages 1--13, 2012.
	
	\bibitem{LM04}
	N.~{Le Bihan} and J.~Mars.
	\newblock Singular value decomposition of quaternion matrices: a new tool for
	vector-sensor signal processing.
	\newblock {\em Signal Processing}, 84:1177--1199, 2004.
	
	\bibitem{LLB13}
	G.~Leclercq, P.~Lef\'{e}vre, and G.~Blohm.
	\newblock 3{D} kinematics using dual quaternions: theory and applications in
	neuroscience.
	\newblock {\em Frontiers in Behavioral Neuroscience}, 7:Article 7 (25pp), 2013.

\bibitem{LHQ22}
	C.~Ling, H.J. He, and L.~Qi.
	\newblock Singular values of dual quaternion matrices and their low-rank approximations.
	\newblock {\em Numerical Functional Analysis and Optimization}, 43 (12):1423--1458, 2022.
	
	\bibitem{LQY22}
	C.~Ling, L.~Qi, and H.~Yan.
	\newblock Minimax principle for right eigenvalues of dual quaternion matrices
	and their generalized inverses.
	\newblock arXiv: 2203.03161, 2022.
	
	\bibitem{MXXLY18}
	Y.~Min, Z.~Xiong, L.~Xing, J.~Liu, and D.~Yin.
	\newblock An improved {SINS/GNSS/CNS} federal filter based on dual quaternions
	(in Chinese).
	\newblock {\em Acta Armamentarii}, 39:315--324, 2018.
	
	\bibitem{Mir59}
	L.~Mirsky.
	\newblock On the trace of matrix products.
	\newblock {\em Mathematische Nachrichten}, 20:171--174, 1959.
	
	\bibitem{Mir75}
	L.~Mirsky.
	\newblock A trace inequality of {John von Neumann}.
	\newblock {\em Monatshefte f\"{u}r Mathematik}, 79:303--306, 1975.
	
	\bibitem{PV10}
	E.~Pennestr\'{i} and P.P. Valentini.
	\newblock Dual quaternions as a tool for rigid body motion analysis: A tutorial
	with an application to biomechanics.
	\newblock {\em The Archive of Mechanical Engineering}, 57:187--205, 2010.
	
	\bibitem{QLY21}
	L.~Qi, C.~Ling, and H.~Yan.
	\newblock Dual quaternions and dual quaternion vectors.
	\newblock {\em Communications on Applied Mathematics and Computation}, 4:1494--1508, 2022.
	
	\bibitem{QL21}
	L.~Qi and Z.~Luo.
	\newblock Eigenvalues and singular values of dual quaternion matrices.
	\newblock {\em Pacific Journal of Optimization}, 19 (2):257--272, 2023.
	
	\bibitem{QWL22}
	L.~Qi, X.~Wang, and Z. Luo.
	\newblock Dual quaternion matrices in multi-agent formation control.
	\newblock to appear in: \newblock {Communications in Mathematical Sciences}.
	
	\bibitem{Tho14}
	F.~Thomas.
	\newblock Approaching dual quaternions from matrix algebra.
	\newblock {\em IEEE Transactions on Robotics}, 30:1037--1048, 2014.
	
	\bibitem{TRA11}
	A.~Torsello, E.~Rodol\`{a}, and A.~Albarelli.
	\newblock Multiview registration via graph diffusion of dual quaternions.
	\newblock In {\em Proceedings of the XXIV IEEE Conference on Computer Vision
		and Pattern Recognition}, pages 2441--2448, 2011.
	
	\bibitem{VN37}
	J.~von Neumann.
	\newblock Some matrix-inequalities and metrization of matric-space.
	\newblock {\em Tomsk Univ. Rev.}, 1:286--299, 1937.
	
	\bibitem{Wa11}
	X.~Wang.
	\newblock {\em Formation Control in Three Dimensional Space and with Nonlinear
		Dynamics}.
	\newblock Ph.D., National University of Defence Technology, Changsha, China,
	2011.
	
	\bibitem{WYL12}
	X.~Wang, C.~Yu, and Z.~Lin.
	\newblock A dual quaternion solution to attitude and position control for rigid
	body coordination.
	\newblock {\em IEEE Transactions on Robotics}, 28:1162--1170, 2012.
	
	\bibitem{WYZ}
	X.~Wang, C.~Yu, and Z.~Zheng.
	\newblock Multiple rigid-bodies rendezvous problem based on unit dual
	quaternions.
	\newblock Manuscript.
	
	\bibitem{WZ14}
	X.~Wang and H.~Zhu.
	\newblock On the comparisons of unit dual quaternion and homogeneous
	transformation matrix.
	\newblock {\em Advances in Applied Clifford Algebras}, 24:213--229, 2014.
	
	\bibitem{WLZZ18}
	M.~Wei, Y.~Li, F.~Zhang, and J.~Zhao.
	\newblock {\em Quaternion Matrix Computations}.
	\newblock NOVA Science Publishers, New York, 2018.
	
	\bibitem{XM15}
	D.~Xu and D.~Mandic.
	\newblock The theory of quaternion matrix derivatives.
	\newblock {\em IEEE Transactions on Signal Processing}, 63:1543--1556, 2015.
	
	\bibitem{Zh97}
	F.~Zhang.
	\newblock Quaternions and matrices of quaternions.
	\newblock {\em Linear Algebra and its Applications}, 251:21--57, 1997.
	
\end{thebibliography}

\end{document}